\documentstyle[amscd,amssymb,verbatim,epsf]{amsart}

\begin{document}

\theoremstyle{plain}
\newtheorem{Thm}{Theorem}
\newtheorem{Cor}{Corollary}
\newtheorem{Con}{Conjecture}
\newtheorem{Main}{Main Theorem}

\newtheorem{Lem}{Lemma}
\newtheorem{Prop}{Proposition}

\theoremstyle{definition}
\newtheorem{Def}{Definition}
\newtheorem{Note}{Note}

\theoremstyle{remark}
\newtheorem{notation}{Notation}
\renewcommand{\thenotation}{}

\errorcontextlines=0
\numberwithin{equation}{section}
\renewcommand{\rm}{\normalshape}%

\title[Ideal-quasi-Cauchy sequences]%
   {Ideal-quasi-Cauchy sequences}

\author{Huseyin Cakalli*, and Bipan Hazarika**\\*Department of Mathematics, Maltepe University, Marmara E\u{g}\.{I}t\.{I}m K\"oy\"u, TR 34857, Maltepe, \.{I}stanbul-Turkey\\Phone:(+90216)6261050 ext:2248, fax:(+90216)6261113\\**Department of Mathematics, Rajiv Gandhi University, Arunachal Pradesh-India}
\address{H\"Usey\.{I}n \c{C}akall\i\\
        Department of Mathematics, Maltepe University, Marmara E\u{g}\.{I}t\.{I}m K\"oy\"u, TR 34857, Maltepe, \.{I}stanbul-Turkey \; \; \; Phone:(+90216)6261050 ext:2248, fax:(+90216)6261113}
\email{hcakalli@@maltepe.edu.tr; hcakalli@@gmail.com}

\address{Bipan Hazarika\\
Department of Mathematics, Rajiv Gandhi University, Rono Hills, Doimukh-791 112, Itanagar, Arunachal Pradesh-India \; \; \; \; \; Phone:(+91360)2278512 , fax:(+91360)2277881}
\email{bh\_rgu@@yahoo.co.in}
\date{\today}

\begin{abstract}
An ideal $I$ is a family of subsets of positive integers $\textbf{N}$ which is closed under taking finite unions and subsets of its elements.
A sequence $(x_n)$ of real numbers is said to be $I$-convergent to a real number $L$, if for each  \;$ \varepsilon> 0$ the set $\{n:|x_{n}-L|\geq \varepsilon\}$ belongs to $I$. We introduce  $I$-ward compactness of a subset of $\textbf{R}$, the set of real numbers, and $I$-ward continuity of a real function in the senses that a subset $E$ of $\textbf{R}$ is $I$-ward compact if any sequence $(x_{n})$ of points in $E$ has an $I$-quasi-Cauchy subsequence, and a real function is $I$-ward continuous if it preserves $I$-quasi-Cauchy sequences where a sequence $(x_{n})$ is called to be $I$-quasi-Cauchy when $(\Delta x_{n})$ is $I$-convergent to $0$. We obtain results related to $I$-ward continuity, $I$-ward compactness, ward continuity, ward compactness, ordinary compactness, ordinary continuity, $\delta$-ward continuity, and slowly oscillating continuity.
\end{abstract}

\maketitle

\section{Introduction}
\normalfont
The concept of continuity and any concept involving continuity play a very important role not only in pure mathematics but also in other branches of sciences involving mathematics especially in computer science, information theory, and biological science.

A subset $E$ of $\textbf{R}$, the set of real numbers, is compact if any open covering of $E$ has a finite subcovering where $\textbf{R}$ is the set of real numbers. This is equivalent to the statement that any sequence $\textbf{x}=(x_{n})$ of points in $E$ has a convergent subsequence whose limit is in $E$. A real function $f$ is continuous if and only if $(f(x_{n}))$ is a convergent sequence whenever $(x_{n})$ is. Regardless of limit, this is equivalent to the statement that  $(f(x_{n}))$ is Cauchy whenever $(x_{n})$ is. Using the idea of continuity of a real function and the idea of compactness in terms of sequences, \c{C}akall\i~ \cite{CakalliForwardcontinuity} introduced the concept of ward (originally the author used the term "forward" there, and later on he preferred using the term "ward" to the term "forward" in \cite{CakalliStatisticalwardcontinuity}, here in this paper we also use the term "ward" in place of forward) continuity in the sense that a function $f$ is ward continuous if it transforms ward convergent to $0$ sequences to ward convergent to $0$ sequences, i.e. $(f(x_{n}))$ is ward convergent to $0$ whenever $(x_{n})$ is ward convergent to $0$, and introduced the concept of ward compactness in the sense that a subset $E$ of $\textbf{R}$ is ward compact if any sequence $\textbf{x}=(x_{n})$ of points in $E$ has a ward convergent to $0$ subsequence $\textbf{z}=(z_{k})=(x_{n_{k}})$ of the sequence $\textbf{x}$ where a sequence  $(z_{k})$ is called ward convergent to $0$ if  $\lim_{k\rightarrow \infty} (z_{k+1}-z_{k})=0$. In \cite{BurtonColeman} Burton and Coleman, and in \cite{CakalliStatisticalquasiCauchysequences} Cakalli used the term "quasi-Cauchy"  in stead of the term "ward convergent to $0$". From now on in this paper we also prefer to using the term "quasi-Cauchy" to using the term "ward convergent to $0$" for convenience. To make no confusion we again note that a sequence of points in $E$ is quasi-Cauchy if $\lim_{k\rightarrow \infty} (z_{k+1}-z_{k})=0$. In terms of quasi-Cauchy we restate the definitions of ward compactness and ward continuity as follows: a function $f$ is ward continuous if it preserves quasi-Cauchy sequences, i.e. $(f(x_{n}))$ is quasi-Cauchy whenever $(x_{n})$ is, and a subset $E$ of $\textbf{R}$ is ward compact if any sequence $\textbf{x}=(x_{n})$ of points in $E$ has a quasi-Cauchy subsequence $\textbf{z}=(z_{k})=(x_{n_{k}})$ of the sequence $\textbf{x}$. In \cite{CakalliDeltaquasiCauchysequences} a real-valued function  defined on a subset $E$ of $\textsf{R}$ is called $\delta$-ward continuous if it preserves $\delta$-quasi Cauchy sequences where a sequence $\textbf{x}=(x_{n})$ is defined to be $\delta$-quasi Cauchy if the sequence $(\Delta x_{n})$ is quasi Cauchy. A subset $E$ of $\textbf{R}$ is said to be $\delta$-ward compact if any sequence $\textbf{x}=(x_{n})$ of points in $E$ has a $\delta$-quasi-Cauchy subsequence $\textbf{z}=(z_{k})=(x_{n_{k}})$ of the sequence $\textbf{x}$.

A sequence $(x_{n})$ of points in $\textbf{R}$ is slowly oscillating if
for any given $\varepsilon>0$, there exists a $\delta=\delta (\varepsilon) >0$ and an $N=N(\varepsilon)$ such that $|x_{m}-x_{n}|<\varepsilon$ if $n\geq N(\varepsilon)$ and $n\leq m \leq (1+\delta)n$ (see \cite{CakalliSlowlyoscillatingcontinuity}). A function defined on a subset $E$ of $\textbf{R}$ is called slowly oscillating continuous if it preserves slowly oscillating sequences, i.e.  $(f(x_n))$ is a slowly oscillating sequence whenever $(x_n)$ is a slowly oscillating sequence of points in $E$. A function defined on a subset $E$ of $\textbf{R}$ is called quasi-slowly oscillating continuous on $E$ if it preserves quasi-slowly oscillating sequences of points in $E$, i.e.  $(f(x_n))$ is a quasi-slowly oscillating sequence whenever $(x_n)$ is a quasi-slowly oscillating sequence of points in $E$ where a sequence $\textbf{x}=(x_{n})$ is called quasi-slowly oscillating if $(\Delta x_{n})=(x_{n+1}-x_{n})$ is a slowly oscillating sequence.

The purpose of this paper is to investigate the concept of ideal ward compactness of a subset of $\textbf{R}$, and the concept of ideal ward continuity of a real function which cannot be given by means of any summability matrix and prove related theorems.

\maketitle

\section{Preliminaries}

First of all, some definitions and notation will be given in the following. Throughout this paper, $\textbf{N}$, and $\textbf{R}$ will denote the set of all positive integers, and the set of all real numbers, respectively. We will use boldface letters $\bf{x},$ $\bf{y},$ $\bf{z},$ ... for sequences $\textbf{x}=(x_{n}),$ $\textbf{y}=(y_{n}),$ $\textbf{z}=(z_{n}),$ ... of terms in $\textbf{R}$. $c$ and $S$ will denote the set of all convergent sequences and the set of all statistically convergent sequences of points in $\textbf{R}$, respectively.

Following the idea given in a 1946 American Mathematical Monthly problem \cite{Buck}, a number of authors Posner \cite{Posner}, Iwinski \cite{Iwinski}, Srinivasan \cite{Srinivasan}, Antoni \cite{Antoni}, Antoni and Salat \cite{AntoniandSalat}, Spigel and Krupnik \cite{SpigelandKrupnik} have studied $A$-continuity defined by a regular summability matrix $A$. Some authors, \"{O}zt\"{u}rk \cite{Ozturk}, Sava\c{s} and Das \cite{SavasandDas}, Borsik and Salat \cite{BorsikSalat} have studied $A$-continuity for methods of almost convergence or for related methods.

The idea of statistical convergence was formerly given under the name "almost convergence" by Zygmund in the first edition of his celebrated monograph published in Warsaw in 1935 \cite{ZygmundTrigonometricseries}. The concept was formally introduced by Fast \cite{Fast} and later was reintroduced by Schoenberg \cite{Schoenberg}, and also independently by Buck \cite{BuckGeneralizedasymptoticdensity}. Although statistical convergence was introduced over nearly the last seventy years, it has become an active area of research for twenty years. This concept has been applied in various areas such as number theory \cite{ErdosSurlesdensitesde}, measure theory \cite{MillerAmeasuretheoresubsequencecharacterizationofstatisticalconvergence}, trigonometric series \cite{ZygmundTrigonometricseries}, summability theory \cite{FreedmanandSemberDensitiesandsummability}, locally convex spaces \cite{MaddoxStatisticalconvergenceinlocallyconvex}, in the study of strong integral summability \cite{ConnorandSwardsonStrongintegralsummabilityandstonecompactification}, turnpike theory \cite{MakarovLevinRubinovMathematicalEconomicTheory}, \cite{MckenzieTurnpiketheory}, \cite{PehlivanandMamedovStatisticalclusterpointsandturnpike}, Banach spaces \cite{ConnorGanichevandKadetsAcharacterizationofBanachspaceswithseparabledualsviaweakstatisticalconvergence}, metrizable topological groups \cite{Cakallilacunarystatisticalconvergenceintopgroups}, and topological spaces \cite{MaioKocinac}, \cite{CakalliandKhan}. It should be also mentioned that the notion of statistical convergence has been considered, in other contexts, by several people like R.A. Bernstein, Z. Frolik, etc. The concept of statistical convergence is a generalization of the usual notion of convergence that, for real-valued sequences, parallels the usual theory of convergence. For a subset $M$ of $\textbf{N}$ the asymptotic density of $M,$ denoted by $\delta(M)$, is given by
\[
\delta(M)=\lim_{n\rightarrow\infty}\frac{1}{n}|\{k\leq n: k\in M\}|,
\]
if this limit exists, where $|\{k \leq n: k\in {M}\}|$ denotes the cardinality of the set $\{k \leq n : k \in{M}\}$. A sequence $(x_{n})$ is statistically convergent to $\ell$ if
\[
\delta(\{n:|x_{n}-\ell|\geq\epsilon\})=0
\]
for every $\epsilon>0$.
In this case $\ell$ is called the statistical limit of $\textbf{x}$. Schoenberg \cite{Schoenberg} studied some basic properties of statistical convergence and also studied the statistical convergence as a summability method. Fridy \cite{FridyOnstatisticalconvergence} gave characterizations of statistical convergence.

By a lacunary sequence $\theta=(k_{r})$, we mean an increasing sequence $\theta=(k_{r})$ of positive integers such that $k_{0}=0$ and $h_{r}:k_{r}-k_{r-1}\rightarrow\infty$. The intervals determined by $\theta$ will be denoted by $I_{r}=(k_{r-1}, k_{r}]$. In this paper, we assume that $\lim inf_{r}\; q_{r}>1$. The notion of lacunary statistical convergence was introduced, and studied by Fridy and Orhan in \cite{FridyandOrhanlacunarystatisconvergence} and \cite{FridyandOrhanlacunarystatisticalsummability} (see also \cite{FreedmanandSemberandRaphaelSomecesarotypesummabilityspaces} and \cite{Cakallilacunarystatisticalconvergenceintopgroups}).
A sequence $(\alpha_{k})$ of points in $\textbf{R}$ is called lacunary statistically convergent to an element $\ell$ of $\textbf{R}$ if
\[
\lim_{r\rightarrow\infty}\frac{1}{h_{r}}|\{k\in I_{r}: |\alpha_{k}-\ell| \geq\varepsilon\}|=0,
\]
for every positive real number $\varepsilon$.
The assumed condition a few lines above ensures the regularity of the lacunary statistical  sequential method.

The notion of $I$-convergence was initially introduced by Kostyrko, $\breve{S}$al$\grave{a}$t and Wilczy$\acute{n}$ski [30]. Although an ideal is defined as a hereditary and additive family of subsets of a non-empty arbitrary set $X$, here in our study it suffices to take $I$ as a family of subsets of $N$, positive integers, i.e. $I\subset2^{\textbf{N}}$, such that  $A\cup B\in I$ for each $ A,B\in I,$ and each subset of an element of $I$ is an element of $I$.
A non-empty family of sets $F\subset2^{\textbf{N}}$ is a filter on $\textbf{N}$ if and only if $\Phi \notin F$, $A\cap B\in F$ for each $A,B\in F,$ and any subset of an element of $F$ is in $F$. An ideal $I$ is called \textit{non-trivial} if $I\neq\Phi$ and $\textbf{N}\notin I.$ Clearly $I$ is a non-trivial ideal if and only if $F=F(I) =\{\textbf{N}-A: A\in I\} $ is a filter in $\textbf{N}$. A non-trivial ideal $I$ is called \textit{admissible} if and only if $\{ \{ n \} :n\in \textbf{N} \}\subset I$. A non-trivial ideal $I$ is maximal if there cannot exists any non-trivial
ideal $J \neq I$ containing $I$ as a subset. Further details on ideals can be found in Kostyrko, et.al (see [30]).
Throughout this paper we assume $I$ is a non-trivial admissible ideal in $\textbf{N}.$ Recall that a sequence $\mathbf{x}=({\it x_{n}})$ of points in $\mathbf{R}$ is said to be $I$-convergent to a real number $\ell$ if $ \{n\in \mathbf{N}: \vert {\it x_n} -\ell\vert \geq\varepsilon \} \in \textit{I}$ for every $\varepsilon >0$. In this case we write $ I-\lim x_n = \ell $. We see that $I$-convergence of a sequence $(x_{n})$ implies $I$-quasi-Cauchyness of $(x_{n})$. Throughout the paper, $I(\textbf{R})$ and $\Delta I$ will denote the set of all $I$-convergent sequences, and the set of all $I$-quasi-Cauchy  sequences of points in $\textbf{R}$, respectively.

If we take $I = I_f= \{A\subseteq \mathbf{N}: A$ is a finite subset of \textbf{N} $\}$, then $I_f$ will be a non-trivial admissible ideal in $\mathbf{N}$ and the corresponding convergence will coincide with the usual convergence. If we take $I = I_\delta = \{A\subseteq \mathbf{N}: \delta(A) =0\}$ where $\delta(A)$ denotes the asymptotic density of the set $ A $, then $I_\delta$ will be a non-trivial admissible ideal of $\mathbf{N}$ and the corresponding convergence will coincide with the statistical convergence.

Connor and Grosse-Erdman \cite{ConnorGrosseErdmann} gave sequential definitions of continuity for real functions calling $G$-continuity instead of $A$-continuity and their results covers the earlier works related to $A$-continuity where a method of sequential convergence, or briefly a method, is a linear function $G$ defined on a linear subspace of $s$, denoted by $c_{G}$, into $\textbf{R}$. A sequence $\textbf{x}=(x_{n})$ is said to be $G$-convergent to $\ell$ if $\textbf{x}\in c_{G}$ and $G(\textbf{x})=\ell$. In particular, $\lim$ denotes the limit function $\lim \textbf{x}=\lim_{n}x_{n}$ on the linear space $c$, and $st-\lim$ denotes the statistical limit function $st-\lim \textbf{x}=st-\lim_{n}x_{n}$ on the linear space $S$. Also $I-\lim$ denotes the $I-$ limit function $I-\lim \textbf{x}=I-\lim_{n}x_{n}$ on the linear space $I(\textbf{R})$. A method $G$ is called regular if every convergent sequence $\textbf{x}=(x_{n})$ is $G$-convergent with $G(\textbf{x})=\lim \textbf{x}$. A method is called subsequential if whenever $\textbf{x}$ is $G$-convergent with $G(\textbf{x})=\ell$, then there is a subsequence $(x_{n_{k}})$ of $\textbf{x}$ with $\lim_{k} x_{n_{k}}=\ell$. Since the ordinary convergence implies ideal convergence, $I$ is a regular sequential method (\cite{Sleziak$I$-continuityintopologicalspaces}) .  Recently, Cakalli  studied new sequential definitions of compactness in \cite{CakalliSequentialdefinitionsofcompactness},\cite{CakalliNewkindsofcontinuities} and slowly oscillating compactness in \cite{CakalliSlowlyoscillatingcontinuity}.

\maketitle
\section{Ideal sequential compactness}

First we recall the definition of $G$-sequentially compactness of a subset $E$ of $\textbf{R}$.
A subset $E$ of $\textbf{R}$ is called $G$-sequentially compact if whenever $(x_n)$ is a sequence of points in $E$ there is subsequence $\textbf{y}=(y_k)=(x_{n_{k}})$ of $(x_n)$ such that $G(\textbf{y})=\lim \textbf{y}$ in $E$ (see \cite{CakalliOnGcontinuity}).
For regular methods any sequentially compact subset $E$ of $\textbf{R}$ is also $G$-sequentially compact and the converse is not always true. For any regular subsequential method $G$, a subset $E$ of $\textbf{R}$ is $G$-sequentially compact if and only if it is sequentially compact in the ordinary sense.

Although  $I$-sequential compactness is a special case of $G$-sequential compactness when $G=lim$, we state the definition of $I$-sequential compactness of a subset $E$ of $\textbf{R}$ as follows:
\begin{Def}
A subset $E$ of $\textbf{R}$ is called $I$-sequentially compact if whenever $(x_n)$ is a sequence of points in $E$ there is an $I$-convergent subsequence $\textbf{y}=(y_k)=(x_{n_{k}})$ of $(x_n)$ such that $I-\lim \textbf{y}$ is in $E.$
\end{Def}

\begin{Lem} \label{LemconvergentsequencesIconvergent}  Sequential method $I$ is regular, i.e.
If $lim x_{n}=\ell$, then $I-lim x_{n}=\ell$.
\end{Lem}
\begin{pf} The proof follows from the fact that $I$ is admissible(See also \cite{KostyrkoSalatWilczynski}).
\end{pf}

\begin{Lem} \label{Lemidealconvergentsequencehasaconvergentsubsequence} Any $I$-convergent sequence of points in $\textbf{R}$ with an $I$-limit $\ell$ has a convergent subsequence with the same limit $\ell$ in the ordinary sense.
\end{Lem}
\begin{pf} See Proposition 3.2. in \cite{SalatTripathyZimon} for a proof.
\end{pf}

We have the following result which states that any nontrivial admissible ideal $I$ is a regular subsequential sequential method.

\begin{Thm} \label{idealmethodisregularandsubsequential}
The sequential method $I$ is regular and subsequential.

\end{Thm}

\begin{pf}   Regularity of $I$ follows from lemma \ref{LemconvergentsequencesIconvergent}, and subsequentiality of $I$ follows from Lemma \ref{Lemidealconvergentsequencehasaconvergentsubsequence}.
\end{pf}

\begin{Thm} \label{Theoeidealsequential compactnesequivalenttoordinarysequentialcompactness}
A subset of $\textbf{R}$ is sequentially compact if and only if it is $I$-sequentially compact.
\end{Thm}

\begin{pf}
The proof easily follows from Corollary 3 on page 597 in \cite{CakalliSequentialdefinitionsofcompactness} and Lemma \ref{Lemidealconvergentsequencehasaconvergentsubsequence}, so is omitted.
\end{pf}

Although $I$-sequential continuity is a special case of $G$-sequential continuity when $G=lim$ (see also Definition 2 in \cite{Sleziak$I$-continuityintopologicalspaces}), we state the definition of $I$-sequential continuity of a function defined on a subset $E$ of $\textbf{R}$ as follows:

\begin{Def}
A function $f: E\rightarrow \textbf{R}$ is $I$-sequentially continuous at a point $x_0$ if, given a sequence $(x_n)$ of points in $E,$ $I-\lim \textbf{x} = x_0$ implies that $I-\lim f(\textbf{x})= f(x_0).$
\end{Def}

\begin{Thm} \label{Theoeidealcontinuityimpliesordinarycontinuity}
Any $I$-sequentially continuous function at a point $x_0$ is continuous at $x_0$ in the ordinary sense.

\end{Thm}

\begin{pf}
Although there is a proof in \cite{Sleziak$I$-continuityintopologicalspaces} we give a different proof for completeness.
Let $f$ be any $I$-sequentially continuous function at point $x_0$. Since any proper admissible ideal is a regular subsequential method by Theorem  \ref{idealmethodisregularandsubsequential}, it follows from Theorem 13 on page 316 in \cite{CakalliOnGcontinuity} that $f$ is continuous in the ordinary sense.

\end{pf}

\begin{Thm} (Theorem 2.2 in \cite{Sleziak$I$-continuityintopologicalspaces}) \label{Theoeordinarycontinuityimpliesidealcontinuity}
Any continuous function at a point $x_0$ is $I$-sequentially continuous at $x_0$.

\end{Thm}

Combining Theorem \ref{Theoeidealcontinuityimpliesordinarycontinuity} and Theorem \ref{Theoeordinarycontinuityimpliesidealcontinuity} we have the following:
\begin{Cor}
A function is $I$-sequentially continuous at  a point $x_0$ if and only if it is continuous at $x_0$.
\end{Cor}

As statistical limit is an $I$-sequential method we get Theorem 2 on page 962 in \cite{CakalliNewkindsofcontinuities}:

\begin{Cor} A function is statistically continuous at  a point $x_0$ if and only if it is continuous at $x_0$ in the ordinary sense.
\end{Cor}

As lacunarily statistical limit is an $I$-sequential method we get Theorem 5 on page 962 in \cite{CakalliNewkindsofcontinuities}:

\begin{Cor} A function is lacunarily statistically continuous at a point $x_0$ if and only if it is continuous at $x_0$ in the ordinary sense.

\end{Cor}

\begin{Cor} For any regular subsequential method $G$ a function is $G$-sequentially continuous at  a point $x_0$, then it is $I$-sequentially continuous at $x_0$ .
\end{Cor}
\begin{pf}
The proof follows from Theorem 13 on page 316 in \cite{CakalliOnGcontinuity}.
\end{pf}

\begin{Cor} Any ward continuous function on a subset $E$ of $\textbf{R}$ is $I$-sequentially continuous on $E$.
\end{Cor}

\begin{Thm} If a function is slowly oscillating continuous on a subset $E$ of $\textbf{R}$, then it is $I$-sequentially continuous on $E$.
\end{Thm}
\begin{pf}
Let $f$  be any slowly oscillating continuous on $E$. It follows from Theorem 2.1 in \cite{CakalliSlowlyoscillatingcontinuity} that $f$ is continuous. By Theorem \ref{Theoeordinarycontinuityimpliesidealcontinuity} we see that $f$ is $I$-sequentially continuous on $E$. This completes the proof.
\end{pf}

\begin{Thm} If a function is $\delta$-ward continuous on a subset $E$ of $\textbf{R}$, then it is $I$-sequentially continuous on $E$.

\end{Thm}
\begin{pf}
Let $f$  be any $\delta$-ward continuous function on $E$. It follows from Corollary 2 on page 399 in \cite{CakalliDeltaquasiCauchysequences} that $f$ is continuous. By Theorem \ref{Theoeordinarycontinuityimpliesidealcontinuity} we obtain that $f$ is $I$-sequentially continuous on $E$. This completes the proof.
\end{pf}

\begin{Cor} If a function is quasi-slowly oscillating continuous on a subset $E$ of $\textbf{R}$, then it is $I$-sequentially continuous on $E$.

\end{Cor}

\begin{pf}
Let $f$  be any quasi-slowly oscillating continuous on $E$. It follows from Theorem 3.2 in \cite{DikandCanak} that $f$ is continuous. By Theorem \ref{Theoeordinarycontinuityimpliesidealcontinuity} we deduce that $f$ is $I$-sequentially continuous on $E$. This completes the proof.
\end{pf}

\maketitle
\section{Ideal ward continuity}

We say that a sequence $\textbf{x}=(x_{n})$ is $I$-ward convergent to a number $\ell$ if $I-\lim_{n\rightarrow \infty} \Delta x_{n}=\ell$ where $\Delta x_{n}=x_{n+1}-x_{n}$. For the special case $\ell=0$ we say that $\textbf{x}$ is ideal quasi-Cauchy, or $I$-quasi-Cauchy, in place of $I$-ward convergent to $0$. Thus a sequence $(x_n)$ of points of $\textbf{R}$ is $I$-quasi-Cauchy if $(\Delta x_n)$ is $I$-convergent to $0$.

Now we give the definitions of $I$-ward compactness of a subset of $\textbf{R}$.

\begin{Def}
A subset $E$ of $\textbf{R}$ is called $I$-ward compact if whenever $\textbf{x}=(x_{n})$ is a sequence of points in $E$ there is an $I$-quasi-Cauchy subsequence of $\textbf{x}$.
\end{Def}

We note that this definition of $I$-ward compactness can not be obtained by any $G$-sequential compactness, i.e. by any summability matrix $A$, even by the summability matrix $A=(a_{nk})$ defined by $a_{nk}=-1$ if $k=n$ and $a_{kn}=1$ if $k=n+1$ and
\[G(x)=I-\lim A\textbf{x}=I-\lim_{k\rightarrow\infty}\sum^{\infty}_{n=1}a_{kn}x_{n}=I-\lim_{k\rightarrow \infty} \Delta x_{k}  \;  \;  \;  \; \; (1)
\] (see \cite{CakalliSequentialdefinitionsofcompactness} for the definition of $G$-sequential compactness). Despite that $G$-sequential compact subsets of $\textbf{R}$ should include the singleton set $\{0\}$, $I$-ward compact subsets of $\textbf{R}$ do not have to include the singleton $\{0\}$.

Firstly, we note that any finite subset of $\textbf{R}$ is $I$-ward compact, union of two $I$-ward compact subsets of $\textbf{R}$ is $I$-ward compact and intersection of any $I$-ward compact subsets of $\textbf{R}$ is $I$-ward compact. Furthermore any subset of an $I$-ward compact set is $I$-ward compact and any bounded subset of $\textbf{R}$ is $I$-ward compact. Any compact subset of $\textbf{R}$ is also $I$-ward compact, and the set $\textbf{N}$ is not $I$-ward compact. We note that any slowly oscillating compact subset of $\textbf{R}$ is $I$-ward compact (see \cite{CakalliSlowlyoscillatingcontinuity}) for the definition of slowly oscillating compactness). These observations suggest to us having the following result.

\begin{Thm} \label{TheowardcompactnescoincideswithIwardcompactness}
A subset $E$ of $\textbf{R}$  is ward compact if and only if it is $I$-ward compact.
\end{Thm}

\begin{pf}

Let us suppose first that $E$ is ward compact. It follows from Lemma 2 on page 1725 in \cite{CakalliStatisticalwardcontinuity} that $E$ is bounded. Then for any sequence $(x_{n})$, there exists a convergent subsequence $(x_{n_{k}})$ of $(x_{n})$ whose limit may be in $E$ or not.  Then the sequence $(\Delta x_{n_{k}})$ is a null sequence. Since $I$ is a regular method, $(\Delta x_{n_{k}})$ is $I$-convergent to $0$, so it is $I$-quasi-Cauchy. Thus $E$ is $I$-ward compact. Now to prove the converse suppose that $E$ is $I$-ward compact. Take any sequence $(x_{n})$ of points in $E$. Then there exists an $I$-quasi-Cauchy subsequence $(x_{n_{k}})$ of $(x_{n})$.  Since $I$ is subsequential there exists a convergent subsequence $(x_{n_{k_{m}}})$ of $(x_{n_{k}})$. The sequence $(x_{n_{k_{m}}})$ is a quasi-Cauchy subsequence of the sequence $(x_{n})$. Thus $E$ is ward compact. This completes the proof of the theorem.
\end{pf}

\begin{Thm}  \label{TheoremboundednesscoincideswithIwardcompactness}
A subset $E$ of $\textbf{R}$ is bounded if and only if it is $I$-ward compact.

\end{Thm}

\begin{pf}
Using an idea in the proof of Lemma 2 on page 1725 in \cite{CakalliStatisticalwardcontinuity} and the  preceding theorem the proof can be obtained
easily so is omitted.
\end{pf}

Now we give the definition of $I$-ward continuity of a real function.

\begin{Def}
A function $f$ is called $I$-ward continuous on a subset $E$ of $\textbf{R}$ if $I-\lim_{n \rightarrow \infty}\Delta f(x_{n}) =0$ whenever $I-\lim_{n \rightarrow \infty}\Delta x_{n} =0$ for a sequence $\textbf{x} = (x_n)$ of terms in $E$.
\end{Def}

We note that this definition of continuity can not be obtained by any $A$-continuity, i.e. by any summability matrix $A$, even by the summability matrix $A=(a_{nk})$ defined by $(1)$ however for this special summability matrix $A$ if $A$-continuity of $f$ at the point $0$ implies $I$-ward continuity of $f$, then $f(0)=0$; and if $I$-ward continuity of $f$ implies $A$-continuity of $f$ at the point $0$, then $f(0)=0$.

We note that sum of two $I$-ward continuous functions is $I$-ward continuous but the product of two $I$-ward continuous functions need not be $I$-ward continuous as it can be seen by considering product of the $I$-ward continuous function $f(x)=x$ with itself.

In connection with $I$-quasi-Cauchy sequences and $I$-convergent sequences the problem arises to investigate the following types of  "continuity" of functions on $\textbf{R}$.

\begin{description}
\item[($\delta i $)] $(x_{n}) \in {\Delta}I \Rightarrow (f(x_{n})) \in {\Delta}I$
\item[($\delta ic$)] $(x_{n}) \in {\Delta}I \Rightarrow (f(x_{n})) \in c$
\item[$(c)$] $(x_{n}) \in c \Rightarrow (f(x_{n})) \in c$
\item[$(c\delta i )$] $(x_{n}) \in c \Rightarrow (f(x_{n})) \in {\Delta}I$
\item[$(i)$] $(x_{n}) \in I(\textbf{R}) \Rightarrow (f(x_{n})) \in I(\textbf{R})$
\end{description}

We see that $(\delta i )$ is $I$-ward continuity of $f,$ $(i)$ is $I$-continuity of $f$ and $(c)$ states the ordinary continuity of $f$. It is easy to see that $(\delta ic)$ implies $(\delta i)$, and $(\delta i )$ does not imply $(\delta i c)$;  and $(\delta i)$ implies $(c\delta i)$, and $(c\delta i)$ does not imply $(\delta i)$; $(\delta ic)$ implies $(c)$ and $(c)$ does not imply $(\delta ic)$; and $(c)$ is equivalent to $(c\delta i)$.

Now we give the implication $(\delta i)$ implies $(i)$, i.e. any $I$-ward continuous function is $I$-sequentially continuous.

\begin{Thm} \label{idealwardcontinuityimpliesidealcontinuity} If $f$ is $I$-ward continuous on a subset $E$ of $\textbf{R}$, then it is $I$-sequentially continuous on $E$.
\end{Thm}
\begin{pf}
Suppose that $f$ is an $I$-ward continuous function on a subset $E$ of $\textbf{R}.$ Let $(x_{n})$ be an $I$-quasi-Cauchy sequence of points in $E$. Then the sequence $$(x_1, x_0, x_2, x_0, x_3, x_0,..., x_{n-1}, x_0, x_n, x_0,...)$$ is an $I$-quasi-Cauchy sequence. Since $f$ is $I$-ward continuous, the sequence $$(y_n)=(f(x_{1}),f(x_{0}),f(x_{2}),f(x_{0}),...,f(x_{n}),f(x_{0}),...)$$ is an $I$-quasi-Cauchy sequence. Therefore $I-\lim_{n\rightarrow \infty} \Delta y_n=0$. Hence  $I-\lim_{n\rightarrow \infty} [f(x_{n})-f(x_{0})]=0$. It follows that the sequence $(f(x_{n}))$ $I$-converges to $f(x_{0})$. This completes the proof of the theorem.

The converse is not always true for the function $f(x)=x^{2}$ is an example since $I-\lim_{n\rightarrow \infty} \Delta x_n =0$ for the sequence $(x_n)=(\sqrt{n}).$ But $I-\lim_{n\rightarrow \infty} \Delta f(x_n) \neq 0,$ because $(f(\sqrt{n}))=(n).$
\end{pf}

\begin{Thm} \label{idealwardcontinuityimpliesordinarycontinuity} If $f$ is $I$-ward continuous on a subset $E$ of $\textbf{R}$, then it is continuous on $E$ in the ordinary sense.
\end{Thm}

\begin{pf} Let $f$ be an $I$-ward continuous function on $E$. By Theorem \ref{idealwardcontinuityimpliesidealcontinuity}, $f$ is $I$-sequentially continuous on $E$. It follows from Theorem \ref{Theoeidealcontinuityimpliesordinarycontinuity} that $f$ is continuous on $E$ in the ordinary sense. Thus the proof is completed.
\end{pf}

\begin{Thm} An $I$-ward continuous image of any $I$-ward compact subset of  $\textbf{R}$ is $I$-ward compact.
\end{Thm}
\begin{pf}
Suppose that $f$ is an $I$-ward continuous function on a subset $E$ of $\textbf{R}$, and $E$ is an $I$-ward compact subset of $\textbf{R}.$ Let $\textbf{y}=(y_n)$ be a sequence of points in $f(E).$ Write $y_{n}=f(x_{n})$ where $x_{n}\in {E}$ for each $n \in{\textbf{N}}$. $I$-ward compactness of $E$ implies that there is an  $I$-quasi-Cauchy subsequence $\textbf{z}=(z_{k})=(x_{n_{k}})$ of $(x_n)$. Write $(t_{k})=(f(z_{k}))$. As $f$ is $I$-ward continuous, $(f(z_k))$ is an $I$-quasi-Cauchy subsequence of $\textbf{y}$. Thus $f(E)$ is $I$-ward compact. This completes the proof of the theorem.
\end{pf}
\begin{Cor} An $I$-ward continuous image of any compact subset of $\textbf{R}$ is compact.
\end{Cor}
\begin{pf}
The proof of this theorem follows from Theorem \ref{Theoeidealcontinuityimpliesordinarycontinuity}.
\end{pf}

\begin{Cor} An $I$-ward continuous image of any bounded subset of $\textbf{R}$ is bounded.
\end{Cor}
\begin{pf}
The proof follows from Theorem \ref{TheoremboundednesscoincideswithIwardcompactness} and Theorem \ref{idealwardcontinuityimpliesordinarycontinuity}.
\end{pf}

\begin{Cor} An $I$-ward continuous image of a $G$-sequentially compact subset of  $\textbf{R}$ is $G$-sequentially compact for any regular subsequential  method $G.$
\end{Cor}

It is well known that any continuous function on a compact subset $E$ of $\textbf{R}$ is uniformly continuous on $E.$ It is also true for a regular subsequential method $G$ that any $I$-ward continuous function on a $G$-sequentially compact subset $E$ of $\textbf{R}$ is also uniformly continuous on $E$ (see \cite{CakalliForwardcontinuity}). Furthermore, for ward continuous functions defined on an $I$-ward compact subset of $\textbf{R}$ we have the following:

\begin{Thm} If $(f_{n})$ is a sequence of $I$-ward continuous functions defined on a subset $E$ of  $\textbf{R}$, and $(f_{n})$ is uniformly convergent to a function $f$, then $f$ is $I$-ward continuous on $E$.
\end{Thm}

\begin{pf}  Let $\varepsilon > 0$ and $(x_n)$ be a sequence of points in $E$ such that $I-\lim_{n\rightarrow \infty} \Delta x_n = 0$. By the uniform convergence of $(f_n)$ there exists a positive integer $N$ such that $|f_{n}(x)-f(x)|<\frac{\varepsilon}{3}$ for all $x \in {E}$ whenever $n\geq N$. By the definition of an ideal, for all $x\in E$, we have
$$\{n\in \textbf{N}: |f_n(x)-f(x)|\geq \frac{\varepsilon}{3}\}\in I.$$

\noindent As $f_N$ is $I$-ward continuous on $E$ we have
$$
\{n\in \textbf{N}: |f_{N}(x_{n+1})-f_{N}(x_{n})|\geq \frac{\varepsilon}{3}\} \in {\it I}.
$$
On the other hand, we have
$$ \{n\in \textbf{N}:|f(x_{n+1})-f(x_{n})| \geq {\varepsilon}\}\subseteq \{\{n \in \textbf{N}: |f(x_{n+1})-f_{N}(x_{n+1})|\geq \frac{\varepsilon}{3}\}$$ $$\cup \{n\in \textbf{N}:|f_{N}(x_{n+1})-f_{N}(x_{n})|\geq \frac{\varepsilon}{3}\} \cup \{n\in \textbf{N}:|f_{N}(x_{n})-f(x_{n})| \geq \frac{\varepsilon}{3}\}\}.\;\;\;\;\;\;(3)$$
Since $I$ is an admissible ideal, so the right hand side of the relation (3) belongs to $I,$ we have
$$ \{n\in \textbf{N}:|f(x_{n+1})-f(x_{n})| \geq {\varepsilon}\}\in \textit{I}.$$
This completes the proof of the theorem.
\end{pf}

\begin{Thm} The set of all $I$-ward continuous functions on a subset $E$ of $\textbf{R}$ is a closed subset of the set of all continuous functions on $E,$ i.e. $\overline{\Delta iwc(E)} = \Delta iwc(E)$ where $\Delta iwc(E)$ is the set of all $I$-ward continuous functions on $E,~ \overline{\Delta iwc(E)}$ denotes the set of all cluster points of $\Delta iwc(E).$
\end{Thm}

\begin{pf}  Let $f$ be an element in $\overline{\Delta iwc(E)}.$ Then there exists sequence $(f_n)$ of points in $\Delta iwc(E)$ such that $\lim_{n\rightarrow \infty} f_n =f.$ To show that $f$ is $I$-ward continuous consider a $I$-quasi-sequence $(x_n)$ of points in $E$.  Since $(f_n)$ converges to $f,$ there exists a positive integer $N$ such that for all $x\in E$ and for all $n\geq N,$~ $|f_{n}(x)-f(x)|<\frac{\varepsilon}{3}.$
By the definition of ideal for all $x\in E$, we have
$$\{n\in \textbf{N}: |f_n(x)-f(x)|\geq \frac{\varepsilon}{3}\}\in I.$$
\noindent As $f_N$ is $I$-ward continuous on $E$ we have
$$
\{n\in \textbf{N}: |f_{N}(x_{n+1})-f_{N}(x_{n})|\geq \frac{\varepsilon}{3}\} \in {\it I}.
$$
On the other hand, we have
$$ \{n\in \textbf{N}:|f(x_{n+1})-f(x_{n})| \geq {\varepsilon}\}\subseteq \{\{n \in \textbf{N}: |f(x_{n+1})-f_{N}(x_{n+1})|\geq \frac{\varepsilon}{3}\}$$ $$\cup \{n\in \textbf{N}:|f_{N}(x_{n+1})-f_{N}(x_{n})|\geq \frac{\varepsilon}{3}\} \cup \{n\in \textbf{N}:|f_{N}(x_{n})-f(x_{n})| \geq \frac{\varepsilon}{3}\}\}.\;\;\;\;\;\;(3)$$
Since $I$ is an admissible ideal, so the right hand side of the relation (3) belongs to $I,$ we have
$$ \{n\in \textbf{N}:|f(x_{n+1})-f(x_{n})| \geq \frac{\varepsilon}{3}\}\in \textit{I}.$$
This completes the proof of the theorem.
\end{pf}

\begin{Cor} The set of all $I$-ward continuous functions on a subset $E$ of $\textbf{R}$ is a complete subspace of the space of all continuous functions on $E.$
\end{Cor}
\begin{pf}
The proof follows from the preceding theorem.
\end{pf}

\section{Conclusion and further problems}

In this paper, two new concepts, namely the concept of $I$-ward continuity of a real function, and the concept of $I$-ward compactness of a subset of \textbf{R} are introduced and investigated. In this investigation we have obtained theorems related to $I$-ward continuity, $I$-ward compactness, compactness, sequential continuity, and uniform continuity. We also introduced and studied some other continuities involving $I$-quasi-Cauchy sequences, statistical sequences, and convergent sequences of points in \textbf{R}. The present work also contains a generalization of results of the paper \cite{CakalliForwardcontinuity}, and some results in \cite{CakalliStatisticalwardcontinuity} and \cite{CakalliStatisticalquasiCauchysequences}.

Finally we note the following further investigation problems arise.

\noindent {\bf Problem 1.} For further study, we suggest to investigate $I$-quasi-Cauchy sequences of fuzzy points and $I$-ward continuity for the fuzzy functions. However due to the change in settings, the definitions and methods of proofs will not always be analogous to these of the present work(for example see \cite{CakalliandPratul}).

\noindent {\bf Problem 2.} For another further study we suggest to introduce a new concept in dynamical sysytems using $I$-ward continuity.

\end{document}